\documentclass[11pt]{article}

\usepackage{latexsym}
\usepackage{amssymb}

\newtheorem{theorem}{Theorem}
\newtheorem{proposition}{Proposition}

\newtheorem{example}{Example}

\begin{document}
{
\begin{center}
{\Large\bf
On the truncated two-dimensional moment problem.}
\end{center}
\begin{center}
{\bf S.M. Zagorodnyuk}
\end{center}
\section{Introduction.}

Moment problems form a classical and fruitful domain of the mathematical analysis. From their origins in 19-th centure 
(the Stieltjes moment problem) there appeared many variations of these problems which became classical: Hamburger, Hausdorff and
trigonometric moment problems,
matrix and operator moment problems, multidimensional moment problems, see classical books~\cite{cit_800_Akhiezer_Krein_Some_questions__Book},
\cite{cit_1000_Shohat_Tamarkin_1943__Book}, \cite{cit_600_Akhiezer_Classical_MP__Book}, \cite{cit_1000_Berezansky_1965__Book}.  

The multidimensional moment problem (both the full and the truncated versions) turned out to be much more complicated than its
one-dimensional prototype~\cite{cit_980_Berg_Christiansen_Ressel__Book}, \cite{cit_990_Marshall_Book}. 
An operator-theoretical interpretation of the multidimensional moment problem was given
by Fuglede in~\cite{cit_980_Fuglede}. It should be noticed that the operator approach to moment problems was introduced
by Naimark in 1940--1943 and then developed by many authors, see historical notes in~\cite{cit_15000_Zagorodnyuk_2017_J_Adv_Math_Stud}.
Elegant conditions for the solvability of the multidimensional moment problem in the case of the support on semi-algebraic sets were given by
Schm\"udgen in~\cite{cit_995_Schmudgen__1991}, \cite{cit_995_Schmudgen__2003}.
Another conditions for the solvability of the multidimensional moment problem, using an extension of the moment sequence, 
were given by Putinar and Vasilescu, see~\cite{cit_993_P_V__1999}, \cite{cit_14100_Vasilescu}.
Developing the idea of Putinar and Vasilescu, we presented another conditions for the solvability of
the two-dimensional moment problem and proposed an algorithm (which essentially consists of solving of linear equations)
for a construction of the solutions set~\cite{cit_14500_Zagorodnyuk_2010_AFA}.
An analytic parametrization for all solutions of the two-dimensional moment problem in a strip was given
in~\cite{cit_14700_Zagorodnyuk_2013_MFAT}. 
Another approach to multidimensional and complex moment problems (including truncated problems), using extension arguments for $*$-semigroups, 
has been developed by Cicho\'n, Stochel and Szafraniec, see~\cite{cit_982_C_St_Sz_2011} and references therein.
Still another approach for the two-dimensional moment problem was proposed by Ovcharenko
in~\cite{cit_991_Ovcharenko_1}, \cite{cit_991_Ovcharenko_2}.

In this paper we shall be focused on the truncated two-dimensional moment problem.
A general approach for this moment problem was given by Curto and Fialkow in their books~\cite{cit_985_Curto_Fialkow__Book1} 
and~\cite{cit_985_Curto_Fialkow__Book2}.
These books entailed a series of papers by a group of mathematicians, see recent papers~\cite{cit_987_Fialkow_2011}, 
\cite{cit_14200_Vasilescu}, \cite{cit_14000_Yoo}
and references therein.
This approach includes an extension of the matrix of prescribed moments with the same rank. While positive extensions are easy to build, the
Hankel-type structure is hard to inherit and this aim needs an involved analysis. 
Effective optimization algorithms for the multidimensional moment problems were given in the book of Lasserre~\cite{cit_985_Lasserre_Book}.
Another approach for truncated moment problems, using a notion of an idempotent, was presented by Vasilescu in~\cite{cit_14220_Vasilescu}.

Consider the following problem:
to find a non-negative measure $\mu(\delta)$, $\delta\in\mathfrak{B}(\mathbb{R}^2)$, such that
\begin{equation}
\label{f1_1}
\int_{\mathbb{R}^2} x_1^m x_2^n d\mu = s_{m,n},\qquad  m\in\mathbb{Z}_{0,M},\quad n\in\mathbb{Z}_{0,N},
\end{equation}
where $\{ s_{m,n} \}_{m\in\mathbb{Z}_{0,M},\ n\in\mathbb{Z}_{0,N}}$ is a prescribed sequence of real numbers; $M,N\in\mathbb{Z}_+$.
This problem is said to be \textit{the truncated two-dimensional moment problem (with rectangular data)}.

Let $K$ be a subset of $\mathbb{R}^2$. The following problem: to find a solution $\mu$ of the truncated two-dimensional moment problem~(\ref{f1_1})
such that
\begin{equation}
\label{f1_2}
\mathrm{supp}\mu \subseteq K,
\end{equation}
is called \textit{the truncated (two-dimensional) $K$-moment problem (with rectangular data)}.
Since no other types of truncations will appear in the sequel, we shall omit
the words about rectangular data.

As a tool for the study of the truncated two-dimensional moment problem we shall use the truncated $K$-moment problem on parallel lines.
Similar to~\cite{cit_14450_Zagorodnyuk_2005}, this allows to consider a set of Hamburger moment problems
and then to analyze the corresponding systems of non-linear inequalities.
For the cases $M=N=1$ and $M=1, N=2$ this approach leads to the necessary and sufficient conditions of the
solvability of the truncated two-dimensional moment problem, while in the cases 
$M=N=2$; $M=2, N=3$; $M=3, N=2$; $M=3, N=3$
some explicit numerical sufficient conditions for the
solvability are obtained. In all these cases a set of solutions (not necessarily atomic) can be
constructed.

\noindent
{\bf Notations. }
As usual, we denote by $\mathbb{R}, \mathbb{C}, \mathbb{N}, \mathbb{Z}, \mathbb{Z}_+$,
the sets of real numbers, complex numbers, positive integers, integers and non-negative integers,
respectively. By $\mathrm{max}\{ a,b \}$ we denote the maximal number of $a$ and $b$.
For arbitrary $k,l\in\mathbb{Z}$ we set 
$$ \mathbb{Z}_{k,l} := \{ j\in\mathbb{Z}:\ k\leq j\leq l \}. $$ 
By $\mathfrak{B}(M)$ we denote the set of all Borel subsets of $M$, where 
$M\subseteq\mathbb{R}$ or $M\subseteq\mathbb{R}^2$.

\section{The truncated two-dimensional moment problems for the cases $M=N=1$ and $M=1, N=2$.}

Choose an arbitrary $N\in\mathbb{Z}_+$ and arbitrary real numbers
$a_j$, $j\in\mathbb{Z}_{0,N}$: $a_0 < a_1 < a_2 < ... < a_N$. Set
\begin{equation}
\label{f2_1}
K_N = K_N(a_0,...,a_N) = \bigcup\limits_{j=0}^N L_j,\qquad L_j := \{ (x_1,x_2)\in\mathbb{R}^2:\ x_2 = a_j \}.
\end{equation}
Thus, $K_N$ is a union on $N+1$ parallel lines in the plane. In this case the $K$-moment problem is reduced to
a set of Hamburger moment problems (cf.~\cite[Theorems 2 and 4]{cit_14450_Zagorodnyuk_2005}).

\begin{proposition}
\label{p2_1}
Let $M,N\in\mathbb{Z}_+$ and $a_j$, $j\in\mathbb{Z}_{0,N}$: $a_0 < a_1 < a_2 < ... < a_N$, be arbitrary.
Consider the truncated $K$-moment problem~(\ref{f1_1}) with $K=K_N(a_0,...,a_N)$.
Let
\begin{equation}
\label{f2_5}
W = W(a_0,a_1,...,a_N) =
\left|
\begin{array}{ccccc} 1 & 1 & ... & 1\\
a_0 & a_1 & ... & a_N\\
a_0^2 & a_1^2 & ... & a_N^2\\
\vdots & \vdots & \ddots & \vdots\\
a_0^N & a_1^N & ... & a_N^N
\end{array}
\right|,
\end{equation}
and $\Delta_{j;m}$ be the determinant obtained from $W$ by replacing $j$-th column with
$$ \left(
\begin{array}{cccc} s_{m,0}\\
s_{m,1}\\
\vdots \\
s_{m,N}\end{array}
\right),\qquad j\in\mathbb{Z}_{0,N},\ m\in\mathbb{Z}_{0,M}. $$
Set
\begin{equation}
\label{f2_7}
s_m(j) := \frac{\Delta_{j;m}}{W},\qquad j\in\mathbb{Z}_{0,N},\ m\in\mathbb{Z}_{0,M}. 
\end{equation}
The truncated $K_N(a_0,a_1,...,a_N)$-moment problem has a solution if and only if
for each $j\in\mathbb{Z}_{0,N}$, the truncated Hamburger moment problem
with moments
$s_m(j)$:
\begin{equation}
\label{f2_9}
\int_\mathbb{R} x^m d\sigma_j = s_m(j),\qquad  m=0,1,...,M,
\end{equation}
is solvable.
Here $\sigma_j$ is a non-negative measure on $\mathfrak{B}(\mathbb{R})$.

Moreover, if $\sigma_j$ is a solution of the Hamburger moment problem~(\ref{f2_9}), $j\in\mathbb{Z}_{0,N}$,
then we may define a measure $\widetilde{\sigma}_j$ by
\begin{equation}
\label{f2_11}
\widetilde{\sigma}_j(\delta) = \sigma_j(\delta \cap \mathbb{R}),\qquad \delta\in\mathfrak{B}(\mathbb{R}^2).
\end{equation}
Here $\mathbb{R}$ means the set $\{ (x_1,x_2)\in\mathbb{R}^2:\ x_2 = 0 \}$.
We define 
\begin{equation}
\label{f2_12}
\widetilde{\sigma_j}'(\delta) = \widetilde{\sigma_j}(\theta_j^{-1}(\delta)),\qquad \delta\in\mathfrak{B}(\mathbb{R}^2),
\end{equation}
where
\begin{equation}
\label{f2_14}
\theta_j((x_1,x_2)) = (x_1,x_2+a_j):\ \mathbb{R}^2 \rightarrow \mathbb{R}^2.
\end{equation}
Then we can define $\mu$ in the following way:
\begin{equation}
\label{f2_15}
\mu(\delta) = \sum_{j=0}^N \widetilde{\sigma_j}'(\delta),\qquad \delta\in\mathfrak{B}(\mathbb{R}^2),
\end{equation}
to get a solution $\mu$ of the truncated $K_N(a_0,a_1,...,a_N)$-moment problem.
\end{proposition}
\noindent
\textbf{Proof.} 
Suppose that the truncated $K_N(a_0,a_1,...,a_N)$-moment problem has a solution $\mu$. 
For an arbitrary $j\in\mathbb{Z}_{0,N}$ we denote: 
$$ \pi_j((x_1,x_2)) = (x_1,x_2-a_j):\ \mathbb{R}^2 \rightarrow \mathbb{R}^2, $$
and
$$ \mu'_j (\delta) = \mu (\pi_j^{-1}(\delta)),\qquad \delta\in\mathfrak{B}(\mathbb{R}^2). $$
Using the measure $\mu'_j (\delta)$ on $\mathfrak{B}(\mathbb{R}^2)$, we define
the measure $\sigma_j$ as a restriction of $\mu'_j (\delta)$ to $\mathfrak{B}(\mathbb{R})$.
Here by $\mathbb{R}$ we mean the set $\{ (x_1,x_2)\in\mathbb{R}^2:\ x_2 = 0 \}$.
With these notations, using the change of variables for measures and the definition of the integral,
for arbitrary
$m\in\mathbb{Z}_{0,M}$, $n\in\mathbb{Z}_{0,N}$, we may write:
$$ s_{m,n} = \int_{\mathbb{R}^2} x_1^m x_2^n d\mu = \sum_{j=0}^N a_j^n \int_{L_j} x_1^m d\mu = 
\sum_{j=0}^N a_j^n \int_\mathbb{R} x_1^m d\mu_j' =
\sum_{j=0}^N a_j^n \int_\mathbb{R} x^m d\sigma_j. $$
Denote $\mathbf{s}_m(j) = \int_\mathbb{R} x^m d\sigma_j$, $j\in\mathbb{Z}_{0,N}$, $m\in\mathbb{Z}_{0,M}$.
Then
\begin{equation}
\label{f2_17}
\left\{ \begin{array}{ccccc}
\mathbf{s}_m(0) + \mathbf{s}_m(1) + \mathbf{s}_m(2) + ... + \mathbf{s}_m(N) = s_{m,0}, \\
a_0 \mathbf{s}_m(0) + a_1 \mathbf{s}_m(1) + a_2 \mathbf{s}_m(2) + ... + a_N \mathbf{s}_m(N) = s_{m,1}, \\           
a_0^2 \mathbf{s}_m(0) + a_1^2 \mathbf{s}_m(1) + a_2^2 \mathbf{s}_m(2) + ... + a_N^2 \mathbf{s}_m(N) = s_{m,2}, \\
\cdots \\
a_0^N \mathbf{s}_m(0) + a_1^N \mathbf{s}_m(1) + a_2^N \mathbf{s}_m(2) + ... + a_N^N \mathbf{s}_m(N) = s_{m,N}, \end{array}
   \right. \quad (m\in\mathbb{Z}_{0,M}).        
\end{equation}
By Cramer's formulas numbers $\mathbf{s}_m(j)$ coincide with numbers $s_m(j)$ from~(\ref{f2_7}). We conclude that
the truncated Hamburger moment problems~(\ref{f2_9}) are solvable.

On the other hand, suppose that the truncated Hamburger moment problems~(\ref{f2_9}) have solutions $\sigma_j$. 
We define measures $\widetilde{\sigma}_j$, $\widetilde{\sigma}_j'$, $\mu$ by~(\ref{f2_11}), (\ref{f2_12}) and (\ref{f2_15}), respectively.
Observe that $\widetilde{\sigma}_j(\mathbb{R}^2\backslash\mathbb{R}) = 0$. Then $\widetilde{\sigma}_j'(\mathbb{R}^2\backslash L_j) = 0$,
and $\mathrm{supp} \mu\subseteq \bigcup_{j=0}^N L_j$.
Using the change of the variable~(\ref{f2_14}) and the definition of $\mu$ we see that
$$ s_m(j) = \int_\mathbb{R} x_1^m d\sigma_j = \int_{L_j} x_1^m d\mu,\qquad j\in\mathbb{Z}_{0,N},\ m\in\mathbb{Z}_{0,M}. $$
Observe that $s_m(j)$ are solutions of the linear system of equations~(\ref{f2_17}). Then
$$ s_{m,n} = \sum_{j=0}^N a_j^n \int_{L_j} x_1^m d\mu =
\int_{\mathbb{R}^2} \sum_{j=0}^N a_j^n \chi_{L_j} (x_1,x_2) x_1^m d\mu = $$
$$ = \int_{\mathbb{R}^2} x_1^m x_2^n d\mu,\qquad m\in\mathbb{Z}_{0,M},\ n\in\mathbb{Z}_{0,N}. $$
Here by $\chi_{L_j}$ we denote the characteristic function of the set $L_j$.
Thus, $\mu$ is a solution of the truncated $K_N(a_0,a_1,...,a_N)$-moment problem.
$\Box$

At first we consider the case $M=1, N=1$ of the truncated two-dimensional moment problem.

\begin{theorem}
\label{t2_1}
Let the truncated two-dimensional moment problem~(\ref{f1_1}) with $M=1,N=1$ and
some $\{ s_{m,n} \}_{m,n\in\mathbb{Z}_{0,1}}$ be given.
This moment problem has a solution if and only if one of the following conditions holds:
\begin{itemize}
\item[$(i)$] $s_{0,0}=s_{0,1}=s_{1,0}=s_{1,1}=0$;

\item[$(ii)$] $s_{0,0}>0$.
\end{itemize}

In the case~$(i)$ the unique solution is $\mu\equiv 0$. In the case~$(ii)$ a solution $\mu$ can be constructed as a solution
of the truncated $K_1(a_0,a_1)$-moment problem with the same $\{ s_{m,n} \}_{m,n\in\mathbb{Z}_{0,1}}$, and
arbitrary $a_0 < \frac{s_{0,1}}{s_{0,0}}$; $a_1 > \frac{s_{0,1}}{s_{0,0}}$.

\end{theorem}
\noindent
\textbf{Proof.} Suppose that the truncated two-dimensional moment problem with $M=N=1$ has a solution $\mu$. 
Of course, $s_{0,0} = \int d\mu \geq 0$.
If
$s_{0,0}=0$ then $\mu\equiv 0$ and condition~$(i)$ holds. 
If $s_{0,0}>0$ then condition~$(ii)$ is true.

On the other hand, if condition~$(i)$ holds then $\mu\equiv 0$ is a solution of the moment problem.
Of course, it is the unique solution (one can repeat the arguments at the beginning of this Proof).
If condition~$(ii)$ holds, choose arbitrary real $a_0,a_1$ such that
$a_0 < \frac{s_{0,1}}{s_{0,0}}$ and $a_1 > \frac{s_{0,1}}{s_{0,0}}$.
Consider the truncated $K_1(a_0,a_1)$-moment problem with $\{ s_{m,n} \}_{m,n\in\mathbb{Z}_{0,1}}$.
Let us check by~Proposition~\ref{p2_1} that this problem is solvable.
We have: $W = a_1 - a_0$,
$$ s_0(0) = \frac{a_1 s_{0,0} - s_{0,1}}{a_1 - a_0} >0,\quad
s_0(1) = \frac{s_{0,1} - a_0 s_{0,0}}{a_1 - a_0} >0, $$
$$ s_1(0) = \frac{a_1 s_{1,0} - s_{1,1}}{a_1 - a_0},\quad
s_1(1) = \frac{s_{1,1} - a_0 s_{1,0}}{a_1 - a_0}. $$
The Hamburger moment problems~(\ref{f2_9}) are solvable~\cite[Theorem 8]{cit_14470_Zagorodnyuk_2010_MFAT}. Their solutions 
can be used to construct a solution~$\mu$ of the truncated two-dimensional moment problem.
$\Box$

We now turn to the case $M=1, N=2$ of the truncated two-dimensional moment problem.

\begin{theorem}
\label{t2_2}
Let the truncated two-dimensional moment problem~(\ref{f1_1}) with $M=1,N=2$ and
some $\{ s_{m,n} \}_{m\in\mathbb{Z}_{0,1},\ n\in\mathbb{Z}_{0,2}}$ be given.
This moment problem has a solution if and only if one of the following conditions holds:
\begin{itemize}
\item[$(a)$] $s_{0,0} = s_{0,1} = s_{0,2} = s_{1,0} = s_{1,1} = s_{1,2} =0$;

\item[$(b)$] $s_{0,0}>0$, and
\begin{equation}
\label{f2_20}
s_{m,n} = \alpha^n s_{m,0},\qquad m=0,1;\ n=1,2,
\end{equation}
for some $\alpha\in\mathbb{R}$.

\item[$(c)$] $s_{0,0}>0$, $s_{0,0} s_{0,2} - s_{0,1}^2 > 0$.

\end{itemize}

In the case~$(a)$ the unique solution is $\mu\equiv 0$. 

In the case~$(b)$ a solution $\mu$ can be constructed as a solution
of the truncated $K_0(\alpha)$-moment problem with moments $\{ s_{m,n} \}_{m\in\mathbb{Z}_{0,1},\ n=0}$.

In the case~$(c)$ a solution $\mu$ can be constructed as a solution
of the truncated $K_2(a_0,a_1,a_2)$-moment problem with the same $\{ s_{m,n} \}_{m\in\mathbb{Z}_{0,1},\ n\in\mathbb{Z}_{0,2}}$,
arbitrary $a_2 > \sqrt{ \frac{s_{0,2}}{s_{0,0}} }$ and $a_1 = \frac{s_{0,1}}{s_{0,0}}$, $a_0 = -a_2$.

\end{theorem}
\noindent
\textbf{Proof.} Suppose that the truncated two-dimensional moment problem with $M=1, N=2$ has a solution $\mu$. 
Choose $p(x_2) = b_0 + b_1 x_2$, where $b_0,b_1$ are arbitrary real numbers. Since
$$ 0 \leq \int p^2 d\mu = s_{0,0} b_0^2 + 2s_{0,1} b_0 b_1 + s_{0,2} b_1^2, $$
then 
$\Gamma_1 := \left(
\begin{array}{cc} s_{0,0} & s_{0,1} \\
s_{0,1} & s_{0,2} \end{array}
\right) \geq 0$.
If $s_{0,0}=0$ then $\mu\equiv 0$ and condition~$(a)$ is true.
If $s_{0,0}>0$ and  
$s_{0,0} s_{0,2} - s_{0,1}^2 =0$, then 
$0$ is an eigenvalue of the matrix $\Gamma_1$
with an eigenvector
$\left(
\begin{array}{cc} c_0 \\
c_1 \end{array}
\right)$, $c_0,c_1\in\mathbb{R}$. Observe that $c_1\not = 0$. Denote $\alpha = -\frac{c_0}{c_1}$.
From the equation $\Gamma_1 \left(
\begin{array}{cc} c_0 \\
c_1 \end{array}
\right) = 0$, it follows that relation~(\ref{f2_20}) holds for $m=0$.
Observe that $\int_{\mathbb{R}^2} (\alpha - x_2)^2 d\mu = 0$. Then
$\mu(\{ (x_1,x_2)\in\mathbb{R}^2:\ x_2\not = \alpha \}) = 0$. For $n=1,2$, we get
$s_{1,n} = \int_{\mathbb{R}^2} x_1 x_2^n d\mu = \alpha^n s_{1,0}$.
Thus, condition~$(b)$ is true. 
Finally, it remains the case~$(c)$.

Conversely, if condition~$(a)$ holds then $\mu\equiv 0$ is a solution of the moment problem.
Since $s_{0,0}=0$, then any solution is equal to $\mu\equiv 0$.

Suppose that condition~$(b)$ holds. Consider the truncated $K_0(\alpha)$-moment problem with moments $\{ s_{m,n} \}_{m\in\mathbb{Z}_{0,1},\ n=0}$.
Let us check by~Proposition~\ref{p2_1} that this problem is solvable.
In fact, $W=1$, $\Delta_{0;m} = s_m(0) = s_{m,0}$, $m=0,1$. 
Since $s_0(0)=s_{0,0} > 0$, then the truncated Hamburger moment problem~(\ref{f2_9}) has a solution.
Then we may construct $\mu$ as it was described in the statement of the theorem.
Remaining moment equalities then follow from relations~(\ref{f2_20}) and the fact that $\mathrm{supp}\mu\subseteq \{ (x_1,x_2)\in\mathbb{R}^2:\ 
x_2 = \alpha \}$.

Suppose that condition~$(c)$ holds. 
Consider the truncated $K_2(a_0,a_1,a_2)$-moment problem with the same $\{ s_{m,n} \}_{m\in\mathbb{Z}_{0,1},\ n\in\mathbb{Z}_{0,2}}$,
arbitrary $a_2 > \sqrt{ \frac{s_{0,2}}{s_{0,0}} }$ and $a_1 = \frac{s_{0,1}}{s_{0,0}}$, $a_0 = -a_2$.
We shall check by~Proposition~\ref{p2_1} that this problem is solvable.
Observe that 
$W(a_0,a_1,a_2) = 2a_2 (a_2^2 - a_1^2) > 0$, and
$$ s_0(0) = \frac{a_2 - a_1}{W}
(a_1 a_2 s_{0,0} - (a_1+a_2) s_{0,1} + s_{0,2}), $$
$$ s_0(1) = \frac{a_2 - a_0}{W}
(-a_0 a_2 s_{0,0} + (a_2+a_0) s_{0,1} + s_{0,2}), $$
$$ s_0(2) = \frac{a_1 - a_0}{W}
(a_0 a_1 s_{0,0} - (a_0+a_1) s_{0,1} + s_{0,2}). $$
For the solvability of the corresponding three truncated Hamburger moment problems it is sufficient
the validity of the following inequalities: $s_0(j)>0$, $j=0,1,2$, which are equivalent to
\begin{equation}
\label{f2_30}
a_1 a_2 s_{0,0} - (a_1 + a_2) s_{0,1} + s_{0,2} > 0,
\end{equation}
\begin{equation}
\label{f2_32}
a_2^2 s_{0,0} - s_{0,2} > 0,
\end{equation}
\begin{equation}
\label{f2_34}
- a_1 a_2 s_{0,0} - (a_1 - a_2) s_{0,1}  + s_{0,2} > 0.
\end{equation}
All these inequalities are true.
Then the solution of the truncated $K_2(a_0,a_1,a_2)$-moment problem exists and provides us with a solution of
the truncated two-dimensional moment problem. $\Box$

\section{The truncated two-dimensional moment problems for the cases $M=N=2$; $M=2, N=3$; $M=3, N=2$; $M=N=3$.}

Consider arbitrary real numbers $\{ s_{m,n} \}_{m,n\in\mathbb{Z}_{0,3} }$, such that
\begin{equation}
\label{f3_100} 
s_{0,0} > 0,\quad s_{0,0} s_{0,2} - s_{0,1}^2 > 0,\quad s_{0,0} s_{2,0} - s_{1,0}^2 > 0.
\end{equation}
Let us study the truncated two-dimensional $K_3(a_0,a_1,a_2,a_3)$-moment problem with 
the moments $\{ s_{m,n} \}_{m,n\in\mathbb{Z}_{0,3} }$ and with
some $a_0 < a_1 < a_2 < a_3$:
\begin{equation}
\label{f3_55} 
a_2\in \left( \frac{ |s_{0,1}| }{ s_{0,0} }, \sqrt{ \frac{ s_{0,2} }{ s_{0,0} } } \right);
\end{equation}
\begin{equation}
\label{f3_57} 
a_3 > \mathrm{max} \left\{ \frac{ \left| s_{0,3} - a_2^2 s_{0,1}  \right| }{ -a_2^2 s_{0,0} + s_{0,2} },
\sqrt{ \frac{ a_2 s_{0,2} + |s_{0,3}|  }{ a_2 s_{0,0} - |s_{0,1}| } } \right\};
\end{equation}
\begin{equation}
\label{f3_05}
a_0 = -a_3,\quad a_1 = -a_2.
\end{equation}
Observe that condition~(\ref{f3_100}) ensures the correctness of all expressions in~(\ref{f3_55}), (\ref{f3_57}).
Let us study by~Proposition~\ref{p2_1}, when this moment problem has a solution.
We have:  
$W = \prod\limits_{1\leq j<i\leq 4} (a_{i-1}-a_{j-1}) > 0$, and for $m\in\mathbb{Z}_{0,3}$,
\begin{equation}
\label{f3_10}
s_m(0) = \frac{ 2a_2 (a_3-a_2)(a_3+a_2) }{W}
\{
- a_2^2 a_3 s_{m,0} + a_2^2 s_{m,1} + a_3 s_{m,2} - s_{m,3}
\},
\end{equation}
\begin{equation}
\label{f3_15}
s_m(1) = - \frac{ (a_2+a_3)(a_3-a_2) 2a_3 }{W}
\{
- a_3^2 a_2 s_{m,0} + a_3^2 s_{m,1} + a_2 s_{m,2} - s_{m,3}
\},
\end{equation}
\begin{equation}
\label{f3_20}
s_m(2) = \frac{ (a_1+a_3)(a_3+a_2) 2a_3 }{W}
\{
a_2 a_3^2 s_{m,0} + a_3^2 s_{m,1} - a_2 s_{m,2} - s_{m,3}
\},
\end{equation}
\begin{equation}
\label{f3_25}
s_m(3) = - \frac{ (a_2+a_3) 2a_2 (a_2+a_3) }{W}
\{
a_3 a_2^2 s_{m,0} + a_2^2 s_{m,1} - a_3 s_{m,2} - s_{m,3}
\}.
\end{equation}
Sufficient conditions for the solvability of the corresponding Hamburger moment problems~(\ref{f2_9}) are 
the following~(\cite[Theorem 8]{cit_14470_Zagorodnyuk_2010_MFAT}):
\begin{equation}
\label{f3_30}
s_0(j) > 0,\quad s_0(j) s_2(j) - (s_1(j))^2 > 0,\qquad j=0,1,2,3. 
\end{equation}
The first inequality in~(\ref{f3_30}) for $j=0,1,2,3$ is equivalent to the following system:
\begin{equation}
\label{f3_35}
\left\{
\begin{array}{cccc}
-a_2^2 a_3 s_{0,0} + a_2^2 s_{0,1} + a_3 s_{0,2} - s_{0,3} > 0 \\
a_3^2 a_2 s_{0,0} - a_3^2 s_{0,1} - a_2 s_{0,2} + s_{0,3} > 0 \\
a_2 a_3^2 s_{0,0} + a_3^2 s_{0,1} - a_2 s_{0,2} - s_{0,3} > 0 \\
- a_3 a_2^2  s_{0,0} - a_2^2 s_{0,1} + a_3 s_{0,2} + s_{0,3} > 0 \end{array}
\right.. 
\end{equation}
The second inequality in~(\ref{f3_30}) for $j=0,1,2,3$ is equivalent to the following inequalities:
$$ ( -a_2^2 a_3 s_{0,0} + a_2^2 s_{0,1} + a_3 s_{0,2} - s_{0,3} )
   ( -a_2^2 a_3 s_{2,0} + a_2^2 s_{2,1} + a_3 s_{2,2} - s_{2,3} ) > $$
$$ > (-a_2^2 a_3 s_{1,0} + a_2^2 s_{1,1} + a_3 s_{1,2} - s_{1,3})^2, $$
$$ ( a_3^2 a_2 s_{0,0} - a_3^2 s_{0,1} - a_2 s_{0,2} + s_{0,3} )
   ( a_3^2 a_2 s_{2,0} - a_3^2 s_{2,1} - a_2 s_{2,2} + s_{2,3} ) > $$
$$ > ( a_3^2 a_2 s_{1,0} - a_3^2 s_{1,1} - a_2 s_{1,2} + s_{1,3})^2, $$
$$ ( a_3^2 a_2 s_{0,0} + a_3^2 s_{0,1} - a_2 s_{0,2} - s_{0,3} )
   ( a_3^2 a_2 s_{2,0} + a_3^2 s_{2,1} - a_2 s_{2,2} - s_{2,3} ) > $$
$$ > (a_3^2 a_2 s_{1,0} + a_3^2 s_{1,1} - a_2 s_{1,2} - s_{1,3})^2, $$
$$ ( -a_2^2 a_3 s_{0,0} - a_2^2 s_{0,1} + a_3 s_{0,2} + s_{0,3} )
   ( -a_2^2 a_3 s_{2,0} - a_2^2 s_{2,1} + a_3 s_{2,2} + s_{2,3} ) > $$
$$ > ( -a_2^2 a_3 s_{1,0} - a_2^2 s_{1,1} + a_3 s_{1,2} + s_{1,3} )^2. $$
Dividing by $a_3$ or $a_3^2$ we obtain that the latter inequalities are equivalent to the following inequalities:
$$ \left( -a_2^2 s_{0,0} + s_{0,2} + \frac{a_2^2 s_{0,1} - s_{0,3}}{a_3}   \right)
   \left( -a_2^2 s_{2,0} + s_{2,2} + \frac{a_2^2 s_{2,1} - s_{2,3}}{a_3} \right) > $$
$$ > \left(-a_2^2 s_{1,0} + s_{1,2} + \frac{a_2^2 s_{1,1} - s_{1,3}}{a_3} \right)^2, $$
$$ \left( a_2 s_{0,0} - s_{0,1} + \frac{-a_2 s_{0,2} + s_{0,3}}{a_3^2} \right)
   \left( a_2 s_{2,0} - s_{2,1} + \frac{-a_2 s_{2,2} + s_{2,3}}{a_3^2} \right) > $$
$$ > \left( a_2 s_{1,0} - s_{1,1} + \frac{-a_2 s_{1,2} + s_{1,3}}{a_3^2} \right)^2, $$
$$ \left( a_2 s_{0,0} + s_{0,1} - \frac{a_2 s_{0,2} + s_{0,3}}{a_3^2} \right)
   \left( a_2 s_{2,0} + s_{2,1} - \frac{a_2 s_{2,2} + s_{2,3}}{a_3^2} \right) > $$
$$ > \left( a_2 s_{1,0} + s_{1,1} - \frac{a_2 s_{1,2} + s_{1,3}}{a_3^2} \right)^2, $$
$$ \left( -a_2^2 s_{0,0} + s_{0,2} + \frac{s_{0,3} - a_2^2 s_{0,1}}{a_3} \right)
   \left( -a_2^2 s_{2,0} + s_{2,2} + \frac{s_{2,3} - a_2^2 s_{2,1}}{a_3} \right) > $$
\begin{equation}
\label{f3_35_1}
> \left( -a_2^2 s_{1,0} + s_{1,2} + \frac{s_{1,3} - a_2^2 s_{1,1}}{a_3} \right)^2. 
\end{equation}
\textit{We additionally assume that}
\begin{equation}
\label{f3_40_1}
(-a_2^2 s_{0,0} + s_{0,2}) (-a_2^2 s_{2,0} + s_{2,2}) > (-a_2^2 s_{1,0} + s_{1,2})^2,
\end{equation}
\begin{equation}
\label{f3_40_2}
(a_2 s_{0,0} - s_{0,1}) (a_2 s_{2,0} - s_{2,1}) > (a_2 s_{1,0} - s_{1,1})^2,
\end{equation}
\begin{equation}
\label{f3_40_3}
(a_2 s_{0,0} + s_{0,1}) (a_2 s_{2,0} + s_{2,1}) > (a_2 s_{1,0} + s_{1,1})^2.
\end{equation}
In this case inequalities~(\ref{f3_35_1}) will be valid, if $a_3$ is sufficiently large.
In fact, inequalities~(\ref{f3_35_1}) have the following obvious structure:
$$ (r_j + \psi_j(a_3)) (l_j + \xi_j(a_3)) > (t_j + \eta_j(a_3))^2,\qquad j\in\mathbb{Z}_{0,3}, $$
while inequalities~(\ref{f3_40_1})-(\ref{f3_40_3}) mean that
$$ r_j l_j > t_j^2,\qquad j\in\mathbb{Z}_{0,3}. $$
Since $\psi_j(a_3)$, $\xi_j(a_3)$ and $\eta_j(a_3)$ tend to zero as $a_3\rightarrow\infty$, then there exists $A = A(a_2) \in\mathbb{R}$
such that inequalities~(\ref{f3_35_1}) hold, if $a_3 > A$.

System~(\ref{f3_35}) can be written in the following form:
\begin{equation}
\label{f3_44}
\left\{
\begin{array}{cc}
\pm ( a_2^2 s_{0,1} - s_{0,3} ) < a_3 ( -a_2^2 s_{0,0} + s_{0,2} ) \\
\pm ( a_3^2 s_{0,1} - s_{0,3} ) < a_2 ( a_3^2 s_{0,0} - s_{0,2} ) \end{array}
\right.. 
\end{equation}
System~(\ref{f3_44}) is equivalent to the following system:
\begin{equation}
\label{f3_45}
\left\{
\begin{array}{cc}
| a_2^2 s_{0,1} - s_{0,3} | < a_3 ( -a_2^2 s_{0,0} + s_{0,2} ) \\
| a_3^2 s_{0,1} - s_{0,3} | < a_2 ( a_3^2 s_{0,0} - s_{0,2} ) \end{array}
\right.. 
\end{equation}
If
\begin{equation}
\label{f3_49}
a_3 > \frac{ | a_2^2 s_{0,1} - s_{0,3} | }{ -a_2^2 s_{0,0} + s_{0,2} },
\end{equation}
and
\begin{equation}
\label{f3_50}
a_3 > \sqrt{ \frac{ | s_{0,3} | + a_2 s_{0,2} }{ a_2 s_{0,0} - | s_{0,1} | } },
\end{equation}
then inequalities~(\ref{f3_45}) are true.
Observe that relation~(\ref{f3_50}) ensures that
$$ a_3^2 |s_{0,1}| + |s_{0,3}| <  a_2 ( a_3^2 s_{0,0} - s_{0,2} ). $$

Quadratic (with respect to $a_3$ or $a_3^2$) inequalities~(\ref{f3_40_1})-(\ref{f3_40_3}) can be verified by elementary means, using
their discriminants. 
Let us apply our considerations to the truncated two-dimensional moment problem.

\begin{theorem}
\label{t3_1}
Let the truncated two-dimensional moment problem~(\ref{f1_1}) with $M=N=3$ and
some $\{ s_{m,n} \}_{m,n\in\mathbb{Z}_{0,3}}$ be given and conditions~(\ref{f3_100}) hold.
Denote by $I_1$, $I_2$ and $I_3$ the sets of positive real numbers $a_2$ satisfying
inequalities~(\ref{f3_40_1}), (\ref{f3_40_2}) and (\ref{f3_40_3}), respectively.
If
\begin{equation}
\label{f3_110}
\left( \frac{ |s_{0,1}| }{ s_{0,0} }, \sqrt{ \frac{ s_{0,2} }{ s_{0,0} } }\right) \cap I_1 \cap I_2 \cap I_3 \not= \emptyset,
\end{equation}
then this moment problem has a solution.

A solution $\mu$ of the moment problem can be constructed as a solution
of the truncated $K_2(-a_3,-a_2,a_2,a_3)$-moment problem with the same $\{ s_{m,n} \}_{m,n\in\mathbb{Z}_{0,3}}$,
with arbitrary $a_2$ from the interval $\left( \frac{ |s_{0,1}| }{ s_{0,0} }, \sqrt{ \frac{ s_{0,2} }{ s_{0,0} } }\right) 
\cap I_1 \cap I_2 \cap I_3$, 
and some positive large $a_3$.

\end{theorem}
\noindent
\textbf{Proof.} 
The proof follows from the preceding considerations.
$\Box$

Let the truncated two-dimensional moment problem~(\ref{f1_1}) with $M,N\in\mathbb{Z}_{2,3}$ and
some $\{ s_{m,n} \}_{m\in\mathbb{Z}_{0,M},\ n\in\mathbb{Z}_{0,N} }$ be given, and conditions~(\ref{f3_100}) hold.
Notice that conditions~(\ref{f3_100}), (\ref{f3_40_1}), (\ref{f3_40_2}), (\ref{f3_40_3})
and the first interval in~(\ref{f3_110}) do not depend on $s_{m,n}$ with indices $m=3$ or $n=3$.
Thus, we can check conditions of Theorem~\ref{t3_1} for this moment problem (keeping undefined moments
as parameters).

\begin{example}
\label{e3_1}
Consider the truncated two-dimensional moment problem~(\ref{f1_1}) with $M=N=2$, and
$$ s_{0,0} = 4ab,\ s_{0,1} = 0,\ s_{0,2} = \frac{4}{3} ab^3,\ s_{1,0} = s_{1,1} = s_{1,2} = 0, $$
$$ s_{2,0} = \frac{4}{3} a^3 b,\ s_{2,1} = 0,\ s_{2,2} = \frac{4}{9} a^3 b^3, $$
where $a,b$ are arbitrary positive numbers.
In this case, condition~(\ref{f3_100}) holds. 
Moreover, we have:
$$ I_1 = (0,+\infty)\backslash\left\{ \frac{1}{\sqrt{3} } b \right\},\quad
I_2 = I_3 = (0,+\infty); $$
$$ \left( \frac{ |s_{0,1}| }{ s_{0,0} }, \sqrt{ \frac{ s_{0,2} }{ s_{0,0} } }\right) =
\left(
0, \frac{1}{\sqrt{3} } b
\right). $$
By~Theorem~\ref{t3_1} we conclude that this moment problem has a solution.

\end{example}

\begin{center}
{\large\bf 
On the truncated two-dimensional moment problem.}
\end{center}
\begin{center}
{\bf S.M. Zagorodnyuk}
\end{center}

We study the truncated two-dimensional moment problem (with rectangular data):
to find a non-negative measure $\mu(\delta)$, $\delta\in\mathfrak{B}(\mathbb{R}^2)$, such that
$\int_{\mathbb{R}^2} x_1^m x_2^n d\mu = s_{m,n}$, $0\leq m\leq M,\quad 0\leq n\leq N$,
where $\{ s_{m,n} \}_{0\leq m\leq M,\ 0\leq n\leq N}$ is a prescribed sequence of real numbers; $M,N\in\mathbb{Z}_+$.
For the cases $M=N=1$ and $M=1, N=2$ explicit numerical necessary and sufficient conditions for the
solvability of the moment problem are given. In the cases 
$M=N=2$; $M=2, N=3$; $M=3, N=2$; $M=3, N=3$
some explicit numerical sufficient conditions for the
solvability are obtained. In all the cases some solutions (not necessarily atomic) of the moment problem can be
constructed.

}

\noindent
Address:

V. N. Karazin Kharkiv National University \newline\indent
School of Mathematics and Computer Sciences \newline\indent
Department of Higher Mathematics and Informatics \newline\indent
Svobody Square 4, 61022, Kharkiv, Ukraine

Sergey.M.Zagorodnyuk@gmail.com; Sergey.M.Zagorodnyuk@univer.kharkov.ua

\end{document}